\documentclass[10pt]{amsart}

\usepackage{amsfonts}
\usepackage{ifthen}
\usepackage{amsthm}
\usepackage{amsmath}
\usepackage{graphicx}
\usepackage{amscd,amssymb,amsthm}
\usepackage{graphicx}
\usepackage{epstopdf}

\newcounter{minutes}
\divide\time by 60
\newcounter{hours}
\multiply\time by 60 \addtocounter{minutes}{-\time}

\newtheorem{lemma}{Lemma}
\newtheorem{theorem}{Theorem}

\newcommand{\real}{\operatorname{Re}}

\title[Radii of starlikeness and convexity of generalized Mittag-Leffler functions]{Radii of starlikeness and convexity of generalized Mittag-Leffler functions}

\author[\'A. Baricz]{\'Arp\'ad Baricz}
\address{Department of Economics, Babe\c{s}-Bolyai University, Cluj-Napoca 400591, Romania}
\address{Institute of Applied Mathematics, \'Obuda University, 1034 Budapest, Hungary}
\email{bariczocsi@yahoo.com}

\author[A. Prajapati]{Anuja Prajapati}
\address{Department of Mathematics, Sambalpur University, Sambalpur 768019, India}
\email{anujaprajapati49@gmail.com}

\begin{document}

\def\thefootnote{}
\footnotetext{ \texttt{File:~\jobname .tex,
          printed: \number\year-0\number\month-0\number\day,
          \thehours.\ifnum\theminutes<10{0}\fi\theminutes}
} \makeatletter\def\thefootnote{\@arabic\c@footnote}\makeatother


\maketitle

\begin{abstract}
In this paper our aim is to find the radii of starlikeness and convexity of the generalized Mittag-Leffler function for three different kinds of normalization by using their Hadamard factorization in such a way that the resulting functions are analytic in the unit disk of the complex plane. The characterization of entire functions from Laguerre-P\'{o}lya class and a result of H. Kumar and M.A. Pathan on the reality of the zeros of generalized Mittag-Leffler functions, which origins goes back to Dzhrbashyan, Ostrovski\u{i} and Peresyolkova, play important roles in this paper. Moreover, the interlacing properties of the zeros of Mittag-Leffler function and its derivative is also useful in the proof of the main results. By using the Euler-Rayleigh inequalities for the real zeros of the generalized Mittag-Leffler function, we obtain some tight lower and upper bounds for the radii of starlikeness and convexity of order zero.
\end{abstract}



\section{Introduction and the Main Results}

The Mittag-Leffler function arises naturally in the solution of fractional order integral equations or fractional order differential equations and especially in the investigations of the fractional generalization of the kinetic equation, random walks, super-diffusive transport and in the study of complex systems. The ordinary and generalized Mittag-Leffler functions interpolate between a purely exponential law and power-law like behavior. The Mittag-Leffler function has been successfully applied in many areas of sciences and engineering. Due to its vast potential of its applications in solving problems on physical, biological, engineering and earth sciences, the Mittag-Leffler function is an important function that finds widespread use in fractional modeling \cite{sergei}. As the exponential function naturally arises in the solution of integer order differential equations, the Mittag-Leffler function plays an analogous role in the solution of non-integer order differential equations. The detailed treatment of Mittag-Leffler functions can be found in \cite{sergei1,gupta,shukla} and in the references therein, see also \cite{ca,fe,ga,le,li} for more details and applications. Geometric properties of Mittag-Leffler functions such as univalence, starlikeness, convexity and close-to-convexity, were established recently in \cite{gonzalez,haubold,bansal}. Special functions have great importance in geometric function theory, and because of this there is an extensive literature dealing with various geometric properties of certain subclasses of analytic univalent functions involving special functions such as Bessel, Struve and Lommel functions of first kind. Baricz and his coauthors investigated in details the determination of the radii of starlikeness and convexity of some normalized forms of these special functions, see for example \cite{baricz, baricz1,baricz2,baricz3,baricz4,baricz5,baricz6,baricz7,baricz8,baricz9} and the references therein for more details. One of the most important things which we have learned in these studies is that the radii of univalence, starlikeness and convexity are obtained as solutions of some transcendental equations and the obtained radii satisfy some interesting inequalities. The positive zeros of special functions and the Laguerre-P\'{o}lya class of real entire functions played an important role in these papers. Motivated by the above series of papers, in this paper our aim is to present some similar results for the generalized three parameter Mittag-Leffler function. For this, three different normalizations are applied in such way that the resulting functions are analytic. By using the Hadamard factorization of the generalized Mittag-Leffler function and combining the methods from \cite{baricz2,baricz1,baricz7,baricz}, we investigate the radii of starlikeness and convexity for each of the three functions. Moreover, we deduce the interlacing properties of the zeros of Mittag-Leffler function and its derivative, and this result is quite useful in the proof of the main results. Furthermore, our aim is also to give some lower and upper bounds for the radii of starlikeness  and convexity of order zero by using Euler-Rayleigh inequalities for the smallest positive zeros of some transcendental equations.

\subsection{Characterization of starlike and convex functions}
In order to present our results we need the following basic definitions. Let $\mathbb{D}_{r}$ be the disk $\{z \in \mathbb{C}:|z|<r\},$ where $r>0.$ Moreover, let $\mathcal{A}$ be the class of analytic functions defined in $\mathbb{D}_r,$ which satisfy the normalization conditions $f(0)=f^{\prime}(0)-1=0.$ Denote by $\mathcal{S}$ the class of functions belonging to $\mathcal{A},$ which are univalent in $\mathbb{D}_{r}.$ The class of starlike functions, denoted by $\mathcal{S}^{*},$ is the subclass of $\mathcal{S},$ which consists of functions $f$ for which the domain $f(\mathbb{D}_r)$ is starlike with respect to $0.$ The analytic description of $\mathcal{S}^{*}$ is
$$\mathcal{S}^{*}=\left\{f \in \mathcal{A}:\real \left(\frac{zf^{\prime}(z)}{f(z)}\right)>0~~ \text{for}~~~ \text{all}~~ z \in \mathbb{D}_r\right\}.$$
Moreover, let $\mathcal{S}^{*}(\rho)$ be the subclass of $\mathcal{S}$ consisting of functions which are starlike of order $\rho$ in $\mathbb{D}_{r},$ where $0\leq \rho <1,$ that is,
 $$\mathcal{S}^{*}(\rho)=\left\{f \in \mathcal{A}:\real \left(\frac{zf^{\prime}(z)}{f(z)}\right)>\rho~~ \text{for}~~~ \text{all}~~ z \in \mathbb{D}_r\right\}.$$
 The real numbers
 $$r^{*}(f)=\sup \left\{r>0:\real \left(\frac{zf^{\prime}(z)}{f(z)}\right)>0~~ \text{for}~~~ \text{all}~~ z \in \mathbb{D}_{r}\right\}$$
 and
 $$r^{*}_{\rho}(f)=\sup \left\{r>0:\real \left(\frac{zf^{\prime}(z)}{f(z)}\right)>\rho~~ \text{for}~~~ \text{all}~~ z \in \mathbb{D}_{r}\right\}$$
 are called the radius of starlikeness and the radius of starlikeness of order $\rho$ of the function $f,$ respectively. We note that $r^{*}(f)=r^{*}_{0}(f)$ is the largest radius such that the image region  $f(\mathbb{D}_{r^{*}(f)})$ is a starlike domain with respect to the origin.

 The class of convex functions, denoted by $\mathcal{C},$ is the subclass of $\mathcal{S}$ which consists of functions $f$ for which the image domain $f(\mathbb{D}_r)$ is a convex domain. The analytic description of $\mathcal{C}$ is
$$\mathcal{C}=\left\{f \in \mathcal{A}:\real \left(1+\frac{zf^{\prime\prime}(z)}{f^{\prime}(z)}\right)>0~~ \text{for}~~~ \text{all}~~ z \in \mathbb{D}_r\right\}.$$
Moreover, let $\mathcal{C}(\rho)$ be the subclass of $\mathcal{S}$ consisting of functions which are convex of order $\rho$ in $\mathbb{D}_{r},$ where $0\leq \rho <1,$ that is,
 $$\mathcal{C}(\rho)=\left\{f \in \mathcal{A}:\real \left(1+\frac{zf^{\prime\prime}(z)}{f^{\prime}(z)}\right)>\rho~~ \text{for}~~~ \text{all}~~ z \in \mathbb{D}_r\right\}.$$
 The real numbers
 $$r^{c}(f)=\sup \left\{r>0:\real \left(1+\frac{zf^{\prime\prime}(z)}{f^{\prime}(z)}\right)>0~~ \text{for}~~~ \text{all}~~ z \in \mathbb{D}_{r}\right\}$$
 and
 $$r^{c}_{\rho}(f)=\sup \left\{r>0:\real \left(1+\frac{zf^{\prime\prime}(z)}{f^{\prime}(z)}\right)>\rho~~ \text{for}~~~ \text{all}~~ z \in \mathbb{D}_{r}\right\}$$
 are called the radius of convexity and the radius of convexity of order $\rho$ of the function $f,$ respectively. We note that $r^{c}(f)=r^{c}_{0}(f)$ is the largest radius such that the image region  $f(\mathbb{D}_{r^{c}(f)})$ is a convex domain.

\subsection{The Laguerre-P\'olya class of real entire functions}
An entire function is a complex-valued function that is holomorphic over the
whole complex plane. An entire function is called real if it maps the real line into itself. A real entire function $\Omega$ belongs to the  Laguerre-P\'{o}lya class $\mathcal{LP}$ if it can be represented in the form $$\Omega(x)=cx^{m}e^{-ax^2+bx}\prod_{n\geq1}\left(1+\frac{x}{x_n}\right)e^{-\frac{x}{x_n}},$$ with $c,b,x_n\in\mathbb{R},$ $a\geq0,$ $m\in\mathbb{N}\cup\{0\}$ and $\sum1/{x_n}^2<\infty.$ We note that the class $\mathcal{LP}$ consists of entire functions which are uniform limits on the compact sets of the complex plane of polynomials with only real zeros. It is important to mention that that $\mathcal{LP}$ is closed under differentiation, that is, if $\Omega\in\mathcal{LP},$ then $\Omega^{(m)}\in\mathcal{LP}$ for each $m\in\{2,3,\dots\}.$

\subsection{The three parameter generalization of the Mittag-Leffler function}
Now, let us consider the function $\phi(\alpha,z)$ defined by
 \begin{equation}
\label{radd1}
\phi(\alpha,z)=\sum_{n\geq 0}\frac{z^{n}}{\Gamma(\alpha n+1)},
 \end{equation}
where $\Gamma$ denotes the Euler gamma function. This function was introduced by Mittag-Leffler in 1903 and therefore it is known as the Mittag-Leffler function.
Another function having similar properties was introduced later by Wiman \cite{wiman} and it is defined by the following series
 \begin{equation}\label{radd2}
 \phi(\alpha,\beta,z)=\sum_{n\geq 0} \frac{z^{n}}{\Gamma(\alpha n+\beta)}.
 \end{equation}
In 1971 Prabhakar \cite{prabhakar} introduced the three parameter function $\phi(\alpha,\beta,\gamma,z)$ in the form of
 \begin{equation}\label{radd3}
 \phi(\alpha,\beta,\gamma,z)=\sum_{n \geq 0}\frac{(\gamma)_{n}z^{n}}{n!\Gamma(\alpha n+\beta)},
 \end{equation}
where $(\gamma)_{n}$ denotes the Pochhammer symbol (or shifted factorial) given in terms of the gamma function by $(a)_n=\Gamma(a+n)/\Gamma(a).$
Some particular cases of $\phi(\alpha,\beta,\gamma,z)$ are
\begin{gather*}
\phi(1,1,1,z)=e^{z}, \qquad
\phi(1,1,2,z)=e^{z}(z+1),\\
\phi(2,1,1,z)=\cosh(\sqrt{z}), \qquad \phi(2,1,2,z)=\cosh(\sqrt{z})+\frac{1}{2}\sqrt{z}\sinh(\sqrt{z}),\\
\phi(2,2,1,z)=\frac{1}{\sqrt{z}}\sinh(\sqrt{z}),\qquad \phi(2,2,2,z)=\frac{1}{2z}(\sqrt{z}\sinh(\sqrt{z})+z\cosh(\sqrt{z})),\\
\phi(2,3,1,z)=\frac{1}{z}\left(\cosh(\sqrt{z})-1 \right),\qquad \phi(2,4,1,z)=\frac{1}{z}\left(\frac{\sinh(\sqrt{z})}{\sqrt{z}}-1 \right).
\end{gather*}

\subsection{A preliminary result on the Mittag-Leffler function}
To state our first preliminary result, we define the following three transformations mapping the set
$$\left\{\left(\frac{1}{\alpha},\beta\right):\alpha>1,\beta>0\right\}$$ into itself:
$$A:\left(\frac{1}{\alpha},\beta\right)\to \left(\frac{1}{2\alpha},\beta\right), \ \ \ B:\left(\frac{1}{\alpha},\beta\right)\to \left(\frac{1}{2\alpha},\alpha+\beta\right),$$
$$C:\left(\frac{1}{\alpha},\beta\right)\to \left\{\begin{array}{ll}\displaystyle\left(\frac{1}{\alpha},\beta-1\right),&\mbox{if}\ \ \beta>1\\ \\ \displaystyle\left(\frac{1}{\alpha},\beta\right),&\mbox{if}\ \ 0<\beta\leq1\end{array}\right..$$
We put $W_b=A(W_a)\cup B(W_a),$ where
$$W_a=\left\{\left(\frac{1}{\alpha},\beta\right):1<\alpha<2,\beta\in[\alpha-1,1]\cup[\alpha,2]\right\},$$
and we denote by $W_i$ the smallest set containing $W_b$ and invariant with respect to $A,$ $B$ and $C,$ that is, if $(a,b)\in W_i,$ then
$A(a,b),$ $B(a,b),$ $C(a,b)\in W_i.$ By using a result of Peresyolkova \cite{peres}, Kumar and Pathan \cite{kumar} recently proved that if
$\left(\frac{1}{\alpha},\beta\right)\in W_i$ and $\gamma>0,$ then all zeros of the generalized Mittag-Leffler function $\phi(\alpha,\beta,\gamma,z)$ are real and negative. It is worth mentioning that
the reality of the zeros as well as their distribution in the case of $\gamma=1,$ that is of Wiman's extension $\phi(\alpha,\beta,z),$ has a rich literature. For more details see for example the papers of Dzhrbashyan \cite{dzhr}, Ostrovski\u{i} and Peresyolkova \cite{ostrov}, Popov and Sedletskii \cite{sedlet}.

The following lemma, which may be of independent interest, plays an important role in the proof of our main results.
\begin{lemma}\label{lemma}
If $\left(\frac{1}{\alpha},\beta\right)\in W_i$ and $\gamma>0,$ then the function $z \mapsto \phi(\alpha,\beta,\gamma,-z^{2})$ has infinitely many zeros which are all real. Denoting by $\lambda_{\alpha,\beta,\gamma,n}$ the $n$th positive zero of $z\mapsto\phi(\alpha,\beta,\gamma,-z^{2}),$ under the same conditions the Weierstrassian decomposition
$$\phi(\alpha,\beta,\gamma,-z^{2})=\frac{1}{\Gamma(\beta)} \prod_{n\geq 1}\left(1-\frac{z^{2}}{\lambda_{\alpha,\beta,\gamma,n}^{2}}\right)$$
is valid. Moreover, if $\xi_{\alpha,\beta,\gamma,n}$ denotes the $n$th positive zero of $\Psi'(\alpha,\beta,\gamma,z),$ where $\Psi(\alpha,\beta,\gamma,z)=z^{\beta}\phi(\alpha,\beta,\gamma,-z^{2}),$ then the positive zeros $\lambda_{\alpha,\beta,\gamma,n}$ and $\xi_{\alpha,\beta,\gamma,n}$ are interlaced.
\end{lemma}

Observe that the function $z \mapsto \phi(\alpha,\beta,\gamma,-z^{2})$ does not belong to $\mathcal{A}.$ Thus first we perform some natural normalization. We define three functions originating from $\phi(\alpha,\beta,\gamma,z)$ as follows:
$$f_{\alpha,\beta,\gamma}(z)=\left(z^{\beta}\Gamma(\beta)\phi(\alpha,\beta,\gamma,-z^{2})\right)^{1/\beta},$$
$$g_{\alpha,\beta,\gamma}(z)=z\Gamma(\beta)\phi(\alpha,\beta,\gamma,-z^{2}),$$
$$h_{\alpha,\beta,\gamma}(z)=z\Gamma(\beta)\phi(\alpha,\beta,\gamma,-z).$$
Obviously these functions belong to the class $\mathcal{A}.$ Of course, there exist infinitely many other normalizations. The main motivation to consider the above ones is that they are similar to the frequently studied normalizations in the literature of Bessel, $q$-Bessel, Struve, Lommel and Wright functions.

\subsection{Radii of starlikeness of the generalized Mittag-Leffler functions}
Now, our aim is to investigate the radii of starlikeness of the normalized forms of the generalized three parameter Mittag-Leffler function, that is of $f_{\alpha,\beta,\gamma},$ $g_{\alpha,\beta,\gamma}$ and $h_{\alpha,\beta,\gamma}.$ Our aim is to show that the radii of starlikeness of order $\rho$ of the generalized three parameter Mittag-Leffler functions are actually solutions of some transcendental equations. Moreover, we will also find some lower and upper bounds for the radii of starlikeness of order zero. Throughout of this paper for simplicity we use the notation $\lambda(\alpha,\beta,\gamma,z)=\phi(\alpha,\beta,\gamma,-z^2).$ The technique of determining the radii of starlikeness in the next theorem follows the ideas from \cite{baricz3}, \cite{baricz1} and \cite{baricz2}. The results of the next theorem are natural extensions of some recent results (see \cite{prajapat}, where the special case of $\gamma=1$ was considered) on the Mittag-Leffler function in \eqref{radd2}.

\begin{theorem}\label{theorem1}
Let $\left(\frac{1}{\alpha},\beta\right)\in W_i,$ $\gamma>0$ and $\rho \in [0,1).$
\begin{enumerate}
\item[\bf a.] The radius of starlikeness of order $\rho$ of $f_{\alpha,\beta,\gamma}$ is $r^{*}_{\rho}(f_{\alpha,\beta,\gamma})=x_{\alpha,\beta,\gamma,1},$ where $x_{\alpha,\beta,\gamma,1}$ is the smallest positive zero of the transcendental equation
$$r\lambda'(\alpha,\beta,\gamma,r)-\beta (\rho-1)\lambda(\alpha,\beta,\gamma,r)=0.$$
\item[\bf b.] The radius of starlikeness of order $\rho$ of $g_{\alpha,\beta,\gamma}$ is $r^{*}_{\rho}(g_{\alpha,\beta,\gamma})=y_{\alpha,\beta,\gamma,1},$ where $y_{\alpha,\beta,\gamma,1}$ is the smallest positive zero of the transcendental equation
$$r\lambda'(\alpha,\beta,\gamma,r)- (\rho-1)\lambda(\alpha,\beta,\gamma,r)=0.$$
\item[\bf c.] The radius of starlikeness of order $\rho$ of $h_{\alpha,\beta,\gamma}$ is $r^{*}_{\rho}(h_{\alpha,\beta,\gamma})=z_{\alpha,\beta,\gamma,1},$ where $z_{\alpha,\beta,\gamma,1}$ is the smallest positive zero of the transcendental equation
$$\sqrt{r}\lambda'(\alpha,\beta,\gamma,\sqrt{r})-2(\rho-1)\lambda(\alpha,\beta,\gamma,\sqrt{r})=0.$$
\end{enumerate}
\end{theorem}

The following theorem provides some tight lower and upper bounds for the radii of starlikeness of the functions considered in the above theorem.
The technique used in the proof of this theorem follows the ideas from \cite{baricz4} and \cite{baricz5}, and the main idea is to deduce some Euler-Rayleigh inequalities for the first positive zero of some entire functions, which are connected with the transcendental equations in the above theorem. We mention that it is possible to get more sharp results in the next theorem by using higher order Euler-Rayleigh inequalities for $k\in\{2,3,\dots\},$ however we omitted them due to their complicated form.

\begin{theorem}\label{theorem2}
Let $\left(\frac{1}{\alpha},\beta\right)\in W_i$ and $\gamma>0.$
\begin{enumerate}
\item[\bf a.] The radius of starlikeness $r^{*}(f_{\alpha,\beta,\gamma})$ satisfies the inequalities
$$\frac{\gamma(\beta+2)\Gamma(\beta)}{\beta\Gamma(\alpha+\beta)}-\frac{(\gamma+1)(\beta+4)\Gamma(\alpha+\beta)}{(\beta+2)\Gamma(2\alpha+\beta)}
<{\left(r^{*}(f_{\alpha,\beta,\gamma})\right)^{-2}}<\frac{\gamma(\beta+2)\Gamma(\beta)}{\beta\Gamma(\alpha+\beta)}.$$
\item[\bf b.] The radius of starlikeness $r^{*}(g_{\alpha,\beta,\gamma})$ satisfies the inequalities
$$\frac{3\gamma\Gamma(\beta)}{\Gamma(\alpha+\beta)}-\frac{5(\gamma+1)\Gamma(\alpha+\beta)}{3\Gamma(2\alpha+\beta)}<
\left(r^{*}(g_{\alpha,\beta,\gamma})\right)^{-2}<\frac{3\gamma\Gamma(\beta)}{\Gamma(\alpha+\beta)}.$$
\item[\bf c.] The radius of starlikeness $r^{*}(h_{\alpha,\beta,\gamma})$ satisfies the inequalities
$$\frac{2\gamma\Gamma(\beta)}{\Gamma(\alpha+\beta)}-\frac{3(\gamma+1)\Gamma(\alpha+\beta)}{2\Gamma(2\alpha+\beta)}<
\left(r^{*}(h_{\alpha,\beta,\gamma})\right)^{-1}<\frac{2\gamma\Gamma(\beta)}{\Gamma(\alpha+\beta)}.$$
\end{enumerate}
\end{theorem}

\subsection{Radii of convexity of the generalized Mittag-Leffler functions}
Now, we are going to investigate the radii of convexity of order $\rho$ of the functions $f_{\alpha,\beta,\gamma},$ $g_{\alpha,\beta,\gamma}$ and $h_{\alpha,\beta,\gamma}.$ In addition, we find tight lower and upper bounds for the radii of convexity of order zero for the functions $g_{\alpha,\beta,\gamma}$ and $h_{\alpha,\beta,\gamma}.$ The technique used in the process of finding the radii of convexity in the next theorem is based on the ideas from \cite{baricz7} and \cite{baricz8}.

\begin{theorem}\label{theorem4}
Let $\left(\frac{1}{\alpha},\beta\right)\in W_i,$ $\gamma>0$ and $\rho \in [0,1).$
\begin{enumerate}
\item[\bf a.] The radius of convexity $r_{\rho}^c(f_{\alpha,\beta,\gamma})$ is the smallest  positive root of  the transcendental equation $(rf_{\alpha,\beta,\gamma}(r))'=\rho f'_{\alpha,\beta,\gamma}(r).$
\item[\bf b.] The radius of convexity $r_{\rho}^c(g_{\alpha,\beta,\gamma})$ is the smallest  positive root of  the transcendental equation $(rg_{\alpha,\beta,\gamma}(r))'=\rho g'_{\alpha,\beta,\gamma}(r).$
\item[\bf c.] The radius of convexity $r_{\rho}^c(h_{\alpha,\beta,\gamma})$ is the smallest  positive root of  the transcendental equation $(rh_{\alpha,\beta,\gamma}(r))'=\rho h'_{\alpha,\beta,\gamma}(r).$
\end{enumerate}
\end{theorem}

Finally, we present some lower and upper bounds for the radii of convexity of the functions $g_{\alpha,\beta,\gamma}$ and $h_{\alpha,\beta,\gamma}$ by using the corresponding Euler-Rayleigh inequalities.
\begin{theorem}\label{theorem5}
Let $\left(\frac{1}{\alpha},\beta\right)\in W_i$ and $\gamma>0.$
\begin{enumerate}
\item[\bf a.] The radius of convexity $r^{c}(g_{\alpha,\beta,\gamma})$ satisfies the inequalities
$$\frac{9\gamma\Gamma(\beta)}{\Gamma(\alpha+\beta)}-\frac{25(\gamma+1)\Gamma(\alpha+\beta)}{9\Gamma(2\alpha+\beta)}<
\left(r^{c}(g_{\alpha,\beta,\gamma})\right)^{-2}<\frac{9\gamma\Gamma(\beta)}{\Gamma(\alpha+\beta)}.$$
\item[\bf b.] The radius of convexity $r^{c}(h_{\alpha,\beta,\gamma})$ satisfies the inequalities
$$\frac{4\gamma\Gamma(\beta)}{\Gamma(\alpha+\beta)}-\frac{9(\gamma+1)\Gamma(\alpha+\beta)}{4\Gamma(2\alpha+\beta)}<
\left(r^{c}(h_{\alpha,\beta,\gamma})\right)^{-1}<\frac{4\gamma\Gamma(\beta)}{\Gamma(\alpha+\beta)}.$$
\end{enumerate}
\end{theorem}

\section{Proofs of the main results}
\setcounter{equation}{0}

\begin{proof}[\bf Proof of Lemma \ref{lemma}]
Recall that Kumar and Pathan \cite{kumar} recently proved that if
$\left(\frac{1}{\alpha},\beta\right)\in W_i$ and $\gamma>0,$ then all zeros of the generalized Mittag-Leffler function $z\mapsto\phi(\alpha,\beta,\gamma,z)$ are real and negative. Moreover, it is well known that $z\mapsto\phi(\alpha,\beta,\gamma,z)$ is an entire function of order $1/\alpha$ (see \cite{sergei1}), and this is a non-integer number and lies in $(0,1)$ if $\left(\frac{1}{\alpha},\beta\right)\in W_i.$ It follows that the generalized Mittag-Leffler function has infinitely many zeros, which are all real, and its infinite product exists. Now, from the infinite product representation we get
\begin{equation}\label{deriv1}
\frac{\Psi'(\alpha,\beta,\gamma,z)}{\Psi(\alpha,\beta,\gamma,z)}=\frac{\beta}{z}+\frac{\lambda'(\alpha,\beta,\gamma,z)}{\lambda(\alpha,\beta,\gamma,z)}=\frac{\beta}{z}+\sum_{n \geq 1}\frac{2z}{z^{2}-\lambda^{2}_{\alpha,\beta,\gamma,n}}.
\end{equation}
Differentiating both sides of \eqref{deriv1}, we have
\begin{equation*}
\frac{d}{dz}\left(\frac{\Psi'(\alpha,\beta,\gamma,z)}{\Psi(\alpha,\beta,\gamma,z)}\right)=-\frac{\beta}{z^{2}}-2\sum_{n \geq 1}\frac{z^{2}+\lambda^{2}_{\alpha,\beta,\gamma,n}}{(z^{2}-\lambda^{2}_{\alpha,\beta,\gamma,n})^{2}},\quad z \neq \lambda_{\alpha,\beta,\gamma,n}.
\end{equation*}
The right hand side of the above expression is real and negative for each $z$ real, $\left(\frac{1}{\alpha},\beta\right)\in W_i$ and $\gamma>0.$ Thus, the quotient on the left side of \eqref{deriv1} is a strictly decreasing function from $+\infty$ to $-\infty$ as $z$ increases through real values over the open interval $(\lambda_{\alpha,\beta,\gamma,n}, \lambda_{\alpha,\beta,\gamma,n+1}),$ $n \in \mathbb{N}.$ Hence the function $z\mapsto\Psi'(\alpha,\beta,\gamma,z)$ vanishes just once between two consecutive zeros of the function $z\mapsto\lambda(\alpha,\beta,\gamma,z).$
\end{proof}

\begin{proof}[\bf Proof of Theorem \ref{theorem1}]
We need to show that the inequalities
\begin{equation}\label{rad1}
\real \left(\frac{zf^{\prime}(z)}{f(z)}\right) \geq \rho,\qquad \real \left(\frac{zg^{\prime}(z)}{g(z)}\right) \geq \rho \qquad \mbox{and} \qquad \real \left(\frac{zh^{\prime}(z)}{h(z)}\right) \geq \rho
\end{equation}
are valid for $z \in \mathbb{D}_{r^{*}_{\rho}(f_{\alpha, \beta,\gamma})},$ $z \in \mathbb{D}_{r^{*}_{\rho}(g_{\alpha, \beta,\gamma})}$
and $z \in \mathbb{D}_{r^{*}_{\rho}(h_{\alpha, \beta,\gamma})},$ respectively, and each of the above inequalities does not  hold in any larger disk. Recall that under the corresponding conditions the zeros of the Mittag-Leffler function $\phi(\alpha,\beta,\gamma,z)$ are all real. Thus, according to Lemma \ref{lemma} the Mittag-Leffler function admits the Weierstrassian decomposition of the form
$$\phi(\alpha,\beta,\gamma,-z^{2})=\frac{1}{\Gamma(\beta)} \prod_{n\geq 1}\left(1-\frac{z^{2}}{\lambda^{2}_{\alpha,\beta,\gamma,n}}\right)$$
and this infinite product is uniformly convergent on each compact subset of $\mathbb{C}.$ Denoting, as above, the above expression by $\lambda(\alpha,\beta,\gamma,z),$ and by logarithmic differentiation we get
\begin{equation*}
\frac{\lambda'(\alpha,\beta,\gamma,z)}{\lambda(\alpha,\beta,\gamma,z)}=\sum _{n \geq 1}\frac{-2z}{\lambda^{2}_{\alpha,\beta,\gamma,n}-z^{2}},
\end{equation*}
which in turn implies that
$$\frac{zf^{\prime}_{\alpha,\beta,\gamma}(z)}{f_{\alpha,\beta,\gamma}(z)}=1-\frac{1}{\beta}\sum _{n \geq 1}\frac{2z^{2}}{\lambda^{2}_{\alpha,\beta,\gamma,n}-z^{2}}, \ \ \
\frac{zg^{\prime}_{\alpha,\beta,\gamma}(z)}{g_{\alpha,\beta,\gamma}(z)}=1-\sum _{n \geq 1}\frac{2z^{2}}{\lambda^{2}_{\alpha,\beta,\gamma,n}-z^{2}}$$ and
$$\frac{zh^{\prime}_{\alpha,\beta,\gamma}(z)}{h_{\alpha,\beta,\gamma}(z)}=1-\sum _{n \geq 1}\frac{z}{\lambda^{2}_{\alpha,\beta,\gamma,n}-z}.$$
We know that \cite{baricz2} if $z \in \mathbb{C}$ and $\theta \in \mathbb{R}$ are such that $\theta > |z|,$ then
\begin{equation}\label{cov}
\frac{|z|}{\theta-|z|} \geq \real \left(\frac{z}{\theta-z}\right).
\end{equation}
Thus the inequality
\begin{equation*}
\frac{|z|^{2}}{\lambda^{2}_{\alpha,\beta,\gamma,n}-|z|^{2}} \geq \real \left(\frac{z^{2}}{\lambda^{2}_{\alpha,\beta,\gamma,n}-z^{2}}\right),
\end{equation*}
is valid for every $\left(\frac{1}{\alpha},\beta\right)\in W_i,$ $\gamma>0,$ $n \in \mathbb{N}$ and $|z|< \lambda_{\alpha,\beta,\gamma,1},$ and therefore
under the same conditions we have that
$$
\real \left(\frac{zf^{\prime}_{\alpha,\beta,\gamma}(z)}{f_{\alpha,\beta,\gamma}(z)}\right)=1-\frac{1}{\beta}\real \left(\sum _{n \geq 1}\frac{2z^{2}}{\lambda^{2}_{\alpha,\beta,\gamma,n}-z^{2}}\right)\geq 1-\frac{1}{\beta}\sum _{n \geq 1}\frac{2|z|^{2}}{\lambda^{2}_{\alpha,\beta,\gamma,n}-|z|^{2}}=\frac{|z|f^{\prime}_{\alpha,\beta,\gamma}(|z|)}{f_{\alpha,\beta,\gamma}(|z|)},
$$
$$
\real \left(\frac{zg^{\prime}_{\alpha,\beta,\gamma}(z)}{g_{\alpha,\beta,\gamma}(z)}\right)=1-\real \left(\sum _{n \geq 1}\frac{2z^{2}}{\lambda^{2}_{\alpha,\beta,\gamma,n}-z^{2}}\right)\geq 1-\sum _{n \geq 1}\frac{2|z|^{2}}{\lambda^{2}_{\alpha,\beta,\gamma,n}-|z|^{2}}=\frac{|z|g^{\prime}_{\alpha,\beta,\gamma}(|z|)}{g_{\alpha,\beta,\gamma}(|z|)}
$$
and
$$
\real \left(\frac{zh^{\prime}_{\alpha,\beta,\gamma}(z)}{h_{\alpha,\beta,\gamma}(z)}\right)=1-\real \left(\sum _{n \geq 1}\frac{z}{\lambda^{2}_{\alpha,\beta,\gamma,n}-z}\right)\geq 1-\sum _{n \geq 1}\frac{|z|}{\lambda^{2}_{\alpha,\beta,\gamma,n}-|z|}=\frac{|z|h^{\prime}_{\alpha,\beta,\gamma}(|z|)}{h_{\alpha,\beta,\gamma}(|z|)},
$$
where equalities are attained only when $z=|z|=r.$ The minimum principle for harmonic functions and the previous inequalities imply that the corresponding inequalities in \eqref{rad1} are valid if and only if we have $|z|<x_{\alpha,\beta,\gamma,1},$ $|z|< y_{\alpha,\beta,\gamma,1}$ and $|z|< z_{\alpha,\beta,\gamma,1},$ respectively, where $x_{\alpha,\beta,\gamma,1},$ $y_{\alpha,\beta.\gamma,1}$ and $z_{\alpha,\beta,\gamma,1}$ are the smallest positive roots of the following equations
\begin{equation*}
\frac{rf^{\prime}_{\alpha,\beta,\gamma}(r)}{f_{\alpha,\beta,\gamma}(r)}=\rho, \quad  \frac{rg^{\prime}_{\alpha,\beta,\gamma}(r)}{g_{\alpha,\beta,\gamma}(r)}=\rho \quad \mbox{and} \quad \frac{rh^{\prime}_{\alpha,\beta,\gamma}(r)}{h_{\alpha,\beta,\gamma}(r)}=\rho,
\end{equation*}
which are equivalent to $$r\lambda'(\alpha,\beta,\gamma,r)-\beta(\rho-1)\lambda(\alpha,\beta,\gamma,r)=0, \ \ r\lambda'(\alpha,\beta,\gamma,r)- (\rho-1)\lambda(\alpha,\beta,\gamma,r)=0$$
and
$$\sqrt{r}\lambda'(\alpha,\beta,\gamma,\sqrt{r})-2(\rho-1)\lambda(\alpha,\beta,\gamma,\sqrt{r})=0.$$
\end{proof}

\begin{proof}[\bf Proof of Theorem \ref{theorem2}]
{\bf a.} The radius of starlikeness of the normalized Mittag-Leffler function $f_{\alpha,\beta,\gamma}$ corresponds to the radius of starlikeness of the function $\Psi(\alpha,\beta,\gamma,z)=z^{\beta}\lambda(\alpha,\beta,\gamma,z).$ The infinite series representation of the function $z\mapsto\Psi'(\alpha,\beta,\gamma,z)$ and its derivative are as follows:
\begin{equation}\label{rad4}
\Psi'(\alpha,\beta,\gamma,z)=\sum_{n \geq 0} \frac{(-1)^{n}(\gamma)_{n}(2n+\beta)z^{2n+\beta-1}}{n! \Gamma(\alpha n +\beta)},
\end{equation}
\begin{equation}\label{rad5}
\Psi''(\alpha,\beta,\gamma,z)=\sum_{n \geq 0} \frac{(-1)^{n}(\gamma)_{n}(2n+\beta)(2n+\beta-1)z^{2n+\beta-2}}{n!\Gamma(\alpha n +\beta)}.
\end{equation}
In view of Lemma \ref{lemma}, the function $z\mapsto \Psi(\alpha,\beta,\gamma,z)$ belongs to the Laguerre-P\'olya class $\mathcal{LP}.$ This class of functions is closed under differentiation, and therefore $z\mapsto \Psi'(\alpha,\beta,\gamma,z)$ belongs also to the Laguerre-P\'{o}lya class $\mathcal{LP}.$ Hence the zeros of the function $z\mapsto \Psi'(\alpha,\beta,\gamma,z)$ are all real, and in fact according to Lemma \ref{lemma} they are interlaced with the zeros of $z\mapsto \Psi(\alpha,\beta,\gamma,z).$ Thus, $\Psi'(\alpha,\beta,\gamma,z)$ can be written as
\begin{equation}\label{rad6}
\Psi'(\alpha,\beta,\gamma,z)=\frac{\beta}{\Gamma(\beta)}z^{\beta-1}\prod_{n \geq 1} \left(1-\frac{z^{2}}{\xi^{2}_{\alpha,\beta,\gamma,n}}\right).
\end{equation}
Logarithmic differentiation of both sides of \eqref{rad6} for $|z|<\xi_{\alpha,\beta,\gamma,1}$ gives
\begin{align}\label{rad7}
\frac{z\Psi''(\alpha,\beta,\gamma,z)}{\Psi'(\alpha,\beta,\gamma,z)}-(\beta-1)&=\sum_{n \geq 1} \frac{-2z^{2}}{\xi^{2}_{\alpha,\beta,\gamma,n}-z^{2}}=-2\sum_{n \geq 1} \sum_{k \geq 0} \frac{z^{2k+2}}{\xi^{2k+2}_{\alpha,\beta,\gamma,n}}\\&=-2\sum_{k \geq 0} \sum_{n \geq 1} \frac{z^{2k+2}}{\xi^{2k+2}_{\alpha,\beta,\gamma,n}}=-2\sum_{k \geq 0}\chi_{k+1}z^{2k+2},\nonumber
\end{align}
where $\chi_{k}=\sum_{n \geq 1}\xi^{-2k}_{\alpha,\beta,\gamma,n}.$ On the other hand, by using \eqref{rad4} and \eqref{rad5} we get
\begin{equation}\label{rad8}
\frac{z\Psi''(\alpha,\beta,\gamma,z)}{\Psi'(\alpha,\beta,\gamma,z)}=\left.\sum_{n\geq 0}a_{n}z^{2n}\right/\sum_{n \geq 0}b_{n}z^{2n},
\end{equation}
where $$a_{n}=\frac{(-1)^{n}(\gamma)_{n}(2n+\beta)(2n+\beta-1)}{n!\Gamma(\alpha n+\beta)}\  \ \ \mbox{and}\ \ \ b_{n}=\frac{(-1)^{n}(\gamma)_{n}(2n+\beta)}{n!\Gamma(\alpha n+\beta)}.$$
By comparing the coefficients of \eqref{rad7}  and \eqref{rad8} we have
$$(\beta-1)b_{0}=a_{0},\ \  (\beta-1)b_{1}-2\chi_{1}a_{0}=a_{1},\ \ (\beta-1)b_{2}-2\chi_{1}b_{1}-2\chi_{2}b_{0}=a_{2},$$
which implies that
$$\chi_{1}=\frac{\gamma(\beta+2) \Gamma(\beta)}{\beta \Gamma(\alpha+\beta)},\ \ \chi_{2}=\frac{\gamma^{2}(\beta+2)^{2}\Gamma^{2}(\beta)}{\beta^{2}\Gamma^{2}(\alpha+\beta)}-\frac{\gamma(\gamma+1)(\beta+4)\Gamma(\beta)}{\beta
\Gamma(2\alpha+\beta)}.$$
By using the Euler-Rayleigh inequalities $$\chi_{k}^{-1/k}< \xi^{2}_{\alpha,\beta,\gamma,1}<\frac{\chi_{k}}{\chi_{k+1}}$$ for $k=1$ we have the inequalities of the first part of the theorem.

{\bf b.} If $\rho=0$ in the second part of Theorem \ref{theorem1}, then we have that the radius of starlikeness of order zero of the function $g_{\alpha,\beta,\gamma}$ is the smallest positive root of the equation $(z\lambda(\alpha,\beta,\gamma,z))'=0.$ Therefore, it is of interest to study the first positive zero of
\begin{equation}\label{rad9}
\omega(\alpha,\beta,\gamma,z)=(z\lambda(\alpha,\beta,\gamma,z))'=\sum_{n \geq 0}\frac{(-1)^{n}(2n+1)(\gamma)_{n}z^{2n}}{n!\Gamma(\alpha n+\beta)}.
\end{equation}
We know that the function $z\mapsto \lambda(\alpha,\beta,\gamma,z)$ belongs to the Laguerre-P\'{o}lya class $\mathcal{LP},$ which is closed under differentiation. Therefore, we get that the function $z\mapsto \omega(\alpha,\beta,\gamma,z)$ belongs to the Laguerre-P\'{o}lya class, and hence all its zeros are real. Suppose that $\zeta_{\alpha,\beta,\gamma,n}$ is the $n$th positive zero $z\mapsto \omega(\alpha,\beta,\gamma,z)$. Then the function $z\mapsto\omega(\alpha,\beta,\gamma,z)$ has the following infinite product representation
\begin{equation}\label{rad10}
\omega(\alpha,\beta,\gamma,z)=\frac{1}{\Gamma(\beta)}\prod_{n \geq 1}\left(1-\frac{z^{2}}{\zeta^{2}_{\alpha,\beta,\gamma,n}}\right),
\end{equation}
since its growth order corresponds to the growth order of the generalized Mittag-Leffler function itself. If we take the logarithmic derivative of both sides of \eqref{rad10}, for $|z|<\zeta_{\alpha,\beta,\gamma,1}$ we have
 \begin{equation}\label{rad11}
\frac{\omega'(\alpha,\beta,\gamma,z)}{\omega(\alpha,\beta,\gamma,z)}=\sum_{n\geq1}\frac{-2z}{\zeta_{\alpha,\beta,\gamma,n}^2-z^2}=
\sum_{n\geq1}\sum_{k\geq0}\frac{-2z^{2k+1}}{\zeta_{\alpha,\beta,\gamma,n}^{2k+2}}=\sum_{k\geq0}\sum_{n\geq1}\frac{-2z^{2k+1}}{\zeta_{\alpha,\beta,\gamma,n}^{2k+2}}
=-2 \sum_{k \geq 0}\delta_{k+1}z^{2k+1},
 \end{equation}
 where $\delta_{k}=\sum_{n \geq 1}\zeta^{-2k}_{\alpha,\beta,\gamma,n}.$ Moreover in view of \eqref{rad9}, we have
 \begin{equation}\label{rad12}
 \frac{\omega'(\alpha,\beta,\gamma,z)}{\omega(\alpha,\beta,\gamma,z)}=-2\left.\sum_{n \geq 0}c_{n}z^{2n+1}\right/\sum_{n \geq 0}d_{n}z^{2n},
 \end{equation}
 where $$c_{n}=\frac{(-1)^{n}(2n+3)(\gamma)_{n+1}}{n!\Gamma(\alpha n+\alpha +\beta)} \ \ \ \mbox{and}\ \ \
 d_{n}=\frac{(-1)^{n}(\gamma)_{n}(2n+1)}{n! \Gamma(\alpha n+\beta)}.$$
Comparing the coefficients in \eqref{rad11} and \eqref{rad12} we have that $\delta_{1}d_{0}=c_{0}$ and
$\delta_{2}d_{0}+\delta_{1}d_{1}=c_{1},$ which yields the following Rayleigh sums
$$\delta_{1}=\frac{3\gamma\Gamma(\beta)}{\Gamma(\alpha+\beta)} \ \ \ \mbox{and}\ \ \ \delta_{2}=\frac{9\gamma^2 \Gamma^{2}(\beta)}{\Gamma^{2}(\alpha+\beta)}-\frac{5\gamma(\gamma+1)\Gamma(\beta)}{\Gamma(2\alpha+\beta)}.$$
By using the Euler-Rayleigh inequalities $$\delta_{k}^{-1/k}< \zeta^{2}_{\alpha,\beta,\gamma,1}<\frac{\delta_{k}}{\delta_{k+1}}$$
for $k=1$ we obtain the inequalities of the second part of the theorem.

{\bf c.} Taking $\delta=0$ in the third part of Theorem \ref{theorem1} we obtain that the radius of starlikeness of order zero of the function $h_{\alpha,\beta,\gamma}$ is the smallest positive root of the equation $(z\lambda(\alpha,\beta,\gamma,\sqrt{z}))^{\prime}=0.$ Therefore, it is of interest to study the first positive zero of
\begin{equation}\label{rad13}
\sigma(\alpha,\beta,\gamma,z)=(z\lambda(\alpha,\beta,\gamma,\sqrt{z}))^{\prime}=\sum_{n \geq 0}\frac{(-1)^{n}(\gamma)_{n}(n+1)z^{n}}{ n!\Gamma(\alpha n+\beta)}.
\end{equation}
We know that the function $z\mapsto \lambda(\alpha,\beta,\gamma,z)$ belongs to the Laguerre-P\'{o}lya class $\mathcal{LP},$ and consequently we get that the function $z\mapsto\sigma(\alpha,\beta,\gamma,z)$ belongs also to the Laguerre-P\'{o}lya class. Hence the zeros of the function $z\mapsto \sigma(\alpha,\beta,\gamma,z)$ are all real. Let $\eta_{\alpha,\beta,\gamma,n}$ be the $n$th positive zero of the function $z\mapsto \sigma(\alpha,\beta,\gamma,z)$. Then the next infinite product representation is valid
\begin{equation}\label{rad14}
\sigma(\alpha,\beta,\gamma,z)=\frac{1}{\Gamma(\beta)}\prod_{n \geq 1}\left(1-\frac{z}{\eta_{\alpha,\beta
,\gamma,n}}\right).
\end{equation}
This is in agreement with the fact that according to Kumar and Pathan \cite{kumar} if
$\left(\frac{1}{\alpha},\beta\right)\in W_i$ and $\gamma>0,$ then all zeros of the generalized Mittag-Leffler function $z\mapsto \phi(\alpha,\beta,\gamma,z)$ are real and negative, and consequently all zeros of $z\mapsto \lambda(\alpha,\beta,\gamma,\sqrt{z})$ and then of $z\mapsto \sigma(\alpha,\beta,\gamma,z)$ are all real and positive.

Logarithmic differentiation of both sides of \eqref{rad14} for $|z|<\eta_{\alpha,\beta,\gamma,1}$ gives
\begin{equation}\label{rad15}
\frac{\sigma'(\alpha,\beta,\gamma,z)}{\sigma(\alpha,\beta,\gamma,z)}=-\sum_{n\geq1}\frac{1}{\eta_{\alpha,\beta,\gamma,n}-z}
=-\sum_{n\geq1}\sum_{k\geq0}\frac{z^k}{\eta^{k+1}_{\alpha,\beta,\gamma,n}}=-\sum_{k\geq0}\sum_{n\geq1}\frac{z^k}{\eta^{k+1}_{\alpha,\beta,\gamma,n}}
=-\sum_{k \geq 0}\theta_{k+1}z^{k},
\end{equation}
where $\theta_{k}=\sum_{n \geq 1}\eta^{-k}_{\alpha,\beta,\gamma,n}.$ On the other hand, logarithmic differentiation of both sides of \eqref{rad13} gives
\begin{equation}\label{rad16}
\frac{\sigma'(\alpha,\beta,\gamma,z)}{\sigma(\alpha,\beta,\gamma,z)}=-\left.\sum_{n \geq 0}u_{n}z^{n}\right/\sum_{n \geq 0}v_{n}z^{n},
\end{equation}
where $$u_{n}=\frac{(-1)^{n}(n+2)(\gamma)_{n+1}}{n!\Gamma(\alpha n+\alpha+\beta)} \ \ \ \mbox{and}\ \ \ v_{n}=\frac{(-1)^{n}(n+1)(\gamma)_{n}}{n!\Gamma(\alpha n+\beta)}.$$
By comparing the coefficients of \eqref{rad15} and \eqref{rad16}, we get the following Rayleigh sums
$$\theta_{1}=\frac{2\gamma\Gamma(\beta)}{\Gamma(\alpha+\beta)}\ \ \ \mbox{and}\ \ \ \theta_{2}=\frac{4\gamma^2\Gamma^{2}(\beta)}{\Gamma^{2}(\alpha+\beta)}-\frac{3\gamma(\gamma+1)\Gamma(\beta)}{\Gamma(2\alpha+\beta)}.$$
By using the Euler-Rayleigh inequalities $$\theta_{k}^{-1/k}< \eta_{\alpha,\beta,\gamma,1}<\frac{\theta_{k}}{\theta_{k+1}}$$ for $k=1$ the corresponding part of this theorem is proved.
\end{proof}

\begin{proof}[\bf Proof of Theorem \ref{theorem4}]
{\bf a.} First note that
 $$1+\frac{zf^{\prime\prime}_{\alpha,\beta,\gamma}(z)}{f^{\prime}_{\alpha,\beta,\gamma}(z)}=1+\frac{z\Psi''(\alpha,\beta,\gamma,z)}{\Psi'(\alpha,\beta,\gamma,z)}+
 \left(\frac{1}{\beta}-1\right)\frac{r\Psi'(\alpha,\beta,\gamma,z)}{\Psi(\alpha,\beta,\gamma,z)},$$
 and let us recall the following infinite product representations from the proof of Theorem \ref{theorem1}
 \begin{gather*}
 \Gamma(\beta)\Psi(\alpha,\beta,\gamma,z)=z^{\beta}\prod_{n \geq 1}\left(1-\frac{z^{2}}{\lambda^{2}_{\alpha,\beta,\gamma,n}}\right),\quad \Gamma(\beta)\Psi'(\alpha,\beta,\gamma,z)=\beta z^{\beta-1}\prod_{n \geq 1}\left(1-\frac{z^{2}}{\xi^2_{\alpha,\beta,\gamma,n}}\right),
 \end{gather*}
 where $\lambda_{\alpha,\beta,\gamma,n}$ and $\xi_{\alpha,\beta,\gamma,n}$ are the $n$th positive roots of $z\mapsto\Psi(\alpha,\beta,\gamma,z)$ and $z\mapsto\Psi'(\alpha,\beta,\gamma,z),$ respectively, as in Lemma \ref{lemma}. Logarithmic differentiation of both sides of the above relations yields
 \begin{gather*}
 \frac{z\Psi'(\alpha,\beta,\gamma,z)}{\Psi(\alpha,\beta,\gamma,z)}=\beta-\sum_{n \geq 1}\frac{2z^{2}}{\lambda^{2}_{\alpha,\beta,\gamma,n}-z^{2}},\quad \frac{z\Psi''(\alpha,\beta,\gamma,z)}{\Psi'(\alpha,\beta,\gamma,z)}=\beta-1-\sum_{n \geq 1}\frac{2z^{2}}{\xi^2_{\alpha,\beta,\gamma,n}-z^{2}},
 \end{gather*}
 which implies that
$$ 1+\frac{zf^{\prime\prime}_{\alpha,\beta,\gamma}(z)}{f^{\prime}_{\alpha,\beta,\gamma}(z)}=1-\left(\frac{1}{\beta}-1\right)\sum_{n \geq 1}\frac{2z^{2}}{\lambda^{2}_{\alpha,\beta,\gamma,n}-z^{2}}-\sum_{n \geq 1}\frac{2z^{2}}{\xi^2_{\alpha,\beta,\gamma,n}-z^{2}}.$$
By using the inequality \eqref{cov} for $\beta \in (0,1]$ we have
\begin{equation}\label{rad17}
\real \left(1+\frac{zf^{\prime\prime}_{\alpha,\beta,\gamma}(z)}{f^{\prime}_{\alpha,\beta,\gamma}(z)}\right)\geq 1-\left(\frac{1}{\beta}-1\right)\sum_{n \geq 1}\frac{2r^{2}}{\lambda^{2}_{\alpha,\beta,\gamma,n}-r^{2}}-\sum_{n \geq 1}\frac{2r^{2}}{\xi^{2}_{\alpha,\beta,\gamma,n}-r^{2}},
\end{equation}
 where $|z|=r.$ Moreover, in view of the following inequality (see \cite[Lemma 2.1]{baricz7})
 $$\alpha \real \left(\frac{z}{a-z}\right)-\real \left(\frac{z}{b-z}\right)\geq \alpha \frac{|z|}{a-|z|}-\frac{|z|}{b-|z|},$$
 where $a>b>0,$ $\alpha\in[0,1],$ $z \in \mathbb{C}$ such that $|z|<b,$ we obtain that \eqref{rad17} is also valid  when $\beta > 1$ for all $z \in \mathbb{D}_{\xi_{\alpha,\beta,\gamma,1}}.$ Here we used that the zeros of $\lambda_{\alpha,\beta,\gamma,n}$ and $\xi_{\alpha,\beta,\gamma,n}$ interlace, according to Lemma \ref{lemma}. Now, the above deduced inequalities imply for $r \in (0,\xi_{\alpha,\beta,\gamma,1})$
 $$\inf_{z \in \mathbb{D}_r}\left\{\real \left(1+\frac{zf^{\prime\prime}_{\alpha,\beta,\gamma}(z)}{f^{\prime}_{\alpha,\beta,\gamma}(z)}\right)\right\}=1+\frac{rf^{\prime\prime}_{\alpha,\beta,\gamma}(r)}{f^{\prime}_{\alpha,\beta,\gamma}(r)}.$$
 The function $u_{\alpha,\beta,\gamma}:(0,\xi_{\alpha,\beta,\gamma,1}) \rightarrow \mathbb{R},$ defined by
 $$u_{\alpha,\beta,\gamma}(r)=1+\frac{rf^{\prime\prime}_{\alpha,\beta,\gamma}(r)}{f^{\prime}_{\alpha,\beta,\gamma}(r)},$$
 is strictly decreasing since
 \begin{align*}
 u^{\prime}_{\alpha,\beta,\gamma}(r)&=-\left(\frac{1}{\beta}-1\right)\sum_{n\geq 1}\frac{4r\lambda^{2}_{\alpha,\beta,\gamma,n}}{(\lambda^{2}_{\alpha,\beta,\gamma,n}-r^{2})^{2}}-
 \sum_{n\geq 1}\frac{4r\xi^{2}_{\alpha,\beta,\gamma,n}}{(\xi^{2}_{\alpha,\beta,\gamma,n}-r^{2})^{2}}\\
 &<\sum_{n\geq 1}\frac{4r\lambda^{2}_{\alpha,\beta,\gamma,n}}{(\lambda^{2}_{\alpha,\beta,\gamma,n}-r^{2})^{2}}-
 \sum_{n\geq 1}\frac{4r\xi^{2}_{\alpha,\beta,\gamma,n}}{(\xi^{2}_{\alpha,\beta,\gamma,n}-r^{2})^{2}}<0
 \end{align*}
 for $r \in (0,\xi_{\alpha,\beta,\gamma,1}),$ where we used again the interlacing property of the zeros stated in Lemma \ref{lemma}. Observe that $\lim_{r\searrow 0} u_{\alpha,\beta,\gamma}(r)=1$ and $\lim_{r\nearrow \xi_{\alpha,\beta,\gamma,1}} u_{\alpha,\beta,\gamma}(r)=-\infty,$ which means that for $z \in \mathbb{D}_{r_{1}}$ we get
 $$\real \left(1+\frac{zf^{\prime\prime}_{\alpha,\beta,\gamma}(z)}{f^{\prime}_{\alpha,\beta,\gamma}(z)}\right)>\rho$$  if and only if $r_{1}$ is the unique root of
 $$1+\frac{zf^{\prime\prime}_{\alpha,\beta,\gamma}(r)}{f^{\prime}_{\alpha,\beta,\gamma}(r)}=\rho$$ situated in $(0,\xi_{\alpha,\beta,\gamma,1}).$

{\bf b.} According to \eqref{rad10} we have
 $$g'_{\alpha,\beta,\gamma}(z)=\prod_{n \geq 1}\left(1-\frac{z^{2}}{\zeta^{2}_{\alpha,\beta,\gamma,n}}\right). $$
 Now, taking logarithmic derivatives on both sides, we get
 $$1+\frac{zg^{\prime\prime}_{\alpha,\beta,\gamma}(z)}{g^{\prime}_{\alpha,\beta,\gamma}(z)}=1-\sum_{n \geq 1}\frac{2z^{2}}{\zeta^{2}_{\alpha,\beta,\gamma,n}-z^{2}}.$$
 Application of the inequality \eqref{cov} implies that $$\real \left(1+\frac{zg^{\prime\prime}_{\alpha,\beta,\gamma}(z)}{g^{\prime}_{\alpha,\beta,\gamma}(z)}\right)\geq 1-\sum_{n \geq 1}\frac{2r^{2}}{\zeta^{2}_{\alpha,\beta,\gamma,n}-r^{2}}, $$
 where $|z|=r.$ Thus, for $r \in (0,\zeta_{\alpha,\beta,\gamma,1}),$ we get
 $$\inf_{z \in \mathbb{D}_r}\left\{\real \left(1+\frac{zg^{\prime\prime}_{\alpha,\beta,\gamma}(z)}{g^{\prime}_{\alpha,\beta,\gamma}(z)}\right)\right\}=
 1+\frac{rg^{\prime\prime}_{\alpha,\beta,\gamma}(r)}{g^{\prime}_{\alpha,\beta,\gamma}(r)}.$$
The function $v_{\alpha,\beta,\gamma}:(0,\zeta_{\alpha,\beta,\gamma,1})\rightarrow \mathbb{R},$ defined by
$$v_{\alpha,\beta,\gamma}(r)=1+\frac{rg^{\prime\prime}_{\alpha,\beta,\gamma}(r)}{g^{\prime}_{\alpha,\beta,\gamma}(r)},$$
is strictly decreasing and take the limits $\lim_{r\searrow 0} v_{\alpha,\beta,\gamma}(r)=1$ and $\lim_{r\nearrow \zeta_{\alpha,\beta,\gamma,1}} v_{\alpha,\beta,\gamma}(r)=-\infty,$ which means that for $z \in \mathbb{D}_{r_{2}}$ we get
 $$\real \left(1+\frac{zg^{\prime\prime}_{\alpha,\beta,\gamma}(z)}{g^{\prime}_{\alpha,\beta,\gamma}(z)}\right)>\rho$$  if and only if $r_{2}$ is the unique root of
 $$1+\frac{zg^{\prime\prime}_{\alpha,\beta,\gamma}(r)}{g^{\prime}_{\alpha,\beta,\gamma}(r)}=\rho$$ situated in $(0,\zeta_{\alpha,\beta,\gamma,1}).$

{\bf c.} According to \eqref{rad14} we have
$$h'_{\alpha,\beta,\gamma}(z)=\prod_{n \geq 1}\left(1-\frac{z}{\eta_{\alpha,\beta,\gamma,n}}\right),$$
which implies that
$$1+\frac{zh^{\prime\prime}_{\alpha,\beta,\gamma}(z)}{h^{\prime}_{\alpha,\beta,\gamma}(z)}=1-\sum_{n \geq 1}\frac{z}{\eta_{\alpha,\beta,\gamma,n}-z}.$$
Let $r \in (0,\eta_{\alpha,\beta,\gamma,1})$ be a fixed number. The minimum principle for harmonic function and inequality \eqref{cov} imply that for $z \in \mathbb{D}_{r}$ we have
\begin{align*}
\real \left(1+\frac{zh^{\prime\prime}_{\alpha,\beta,\gamma}(z)}{h^{\prime}_{\alpha,\beta,\gamma}(z)}\right)
&=\real \left(1-\sum_{n \geq 1}\frac{z}{\eta_{\alpha,\beta,\gamma,n}-z}\right)\geq
\min_{|z|=r}\real \left(1-\sum_{n \geq 1}\frac{z}{\eta_{\alpha,\beta,\gamma,n}-z}\right)\\
&=\min_{|z|=r}\left(1-\sum_{n \geq 1}\real \frac{z}{\eta_{\alpha,\beta,\gamma,n}-z}\right)
\geq 1-\sum_{n \geq 1}\frac{r}{\eta_{\alpha,\beta,\gamma,n}-r}
=1+\frac{rh^{\prime\prime}_{\alpha,\beta,\gamma}(r)}{h^{\prime}_{\alpha,\beta,\gamma}(r)}.
\end{align*}
Consequently, it follows that
$$\inf_{z \in \mathbb{D}_r}\left\{\real \left(1+\frac{zh^{\prime\prime}_{\alpha,\beta,\gamma}(z)}{h^{\prime}_{\alpha,\beta,\gamma}(z)}\right)\right\}=1+\frac{rh^{\prime\prime}_{\alpha,\beta,\gamma}(r)}{h^{\prime}_{\alpha,\beta,\gamma}(r)}.$$
Now, let $r_{3}$ be the smallest positive root of the equation
\begin{equation}\label{rad18}
1+\frac{rh^{\prime\prime}_{\alpha,\beta,\gamma}(r)}{h^{\prime}_{\alpha,\beta,\gamma}(r)}=\rho.
\end{equation}
For $z \in \mathbb{D}_{r_{3}},$ we have
$$\real\left(1+\frac{rh^{\prime\prime}_{\alpha,\beta,\gamma}(r)}{h^{\prime}_{\alpha,\beta,\gamma}(r)}\right)>\rho.$$
In order to complete the proof, we need to show that equation \eqref{rad18} has a unique root in $(0,\eta_{\alpha,\beta,\gamma,1}).$ But the equation \eqref{rad18} is equivalent to
$$w_{\alpha,\beta,\gamma}(r)=1-\rho-\sum_{n\geq 1} \frac{r}{\eta_{\alpha,\beta,\gamma,n}-r}=0$$
and we have
$\lim_{r \searrow 0}w_{\alpha,\beta,\gamma}(r)=1-\rho>0,$ and $\lim_{r \nearrow \eta_{\alpha,\beta,\gamma}}w_{\alpha,\beta,\gamma}(r)=-\infty.$
Now, since the function $w_{\alpha,\beta,\gamma}:(0,\eta_{\alpha,\beta,\gamma,1})\to\mathbb{R},$ defined above, is strictly decreasing, it follows that the equation $w_{\alpha,\beta,\gamma}(r)=0$
has indeed a unique root in $(0,\eta_{\alpha,\beta,\gamma,1}).$ This completes the proof of the theorem.
\end{proof}

\begin{proof}[\bf Proof of Theorem \ref{theorem5}]
{\bf a.} By using the infinite series representations of the generalized Mittag-Leffler function and its derivative we obtain
$$\varphi(\alpha,\beta,\gamma,z)=(zg^{\prime}_{\alpha,\beta,\gamma}(z))^{\prime}=1+\sum_{n \geq 1}\frac{(-1)^{n}(\gamma)_{n}\Gamma(\beta)(2n+1)^{2}z^{2n}}{n!\Gamma(\alpha n+\beta)}.$$
We know that $g_{\alpha,\beta,\gamma}\in\mathcal{LP}$ and this in turn implies that $z\mapsto\varphi(\alpha,\beta,\gamma,z)$ belongs also to the Laguerre-P\'{o}lya class and consequently all its zeros are real. Assume that $\tau_{\alpha,\beta,\gamma,n}$ is the $n$th positive zero of the function $z\mapsto\varphi(\alpha,\beta,\gamma,z).$ Then we have that
$$\varphi(\alpha,\beta,\gamma,z)=\prod_{n \geq 1}\left(1-\frac{z^{2}}{\tau^{2}_{\alpha,\beta,\gamma,n}}\right)$$
and for $|z|<\tau_{\alpha,\beta,\gamma,1}$
\begin{equation}\label{rad19}
\frac{\varphi'(\alpha,\beta,\gamma,z)}{\varphi(\alpha,\beta,\gamma,z)}=\sum_{n\geq1}\frac{-2z}{\tau_{\alpha,\beta,\gamma,n}^2-z^2}=
\sum_{n\geq1}\sum_{k\geq0}\frac{-2z^{2k+1}}{\tau_{\alpha,\beta,\gamma,n}^{2k+2}}=\sum_{k\geq0}\sum_{n\geq1}\frac{-2z^{2k+1}}{\tau_{\alpha,\beta,\gamma,n}^{2k+2}}
=-2 \sum_{k \geq 0}\mu_{k+1}z^{2k+1},
\end{equation}
where $\mu_{k}=\sum_{n \geq 1}\tau^{-2k}_{\alpha,\beta,\gamma,n}.$ On the other hand, we have
 \begin{equation}\label{rad20}
 \frac{\varphi'(\alpha,\beta,\gamma,z)}{\varphi(\alpha,\beta,\gamma,z)}=-2\left. \sum_{n \geq 0}q_{n}z^{2n+1}\right/\sum_{n \geq 0} r_{n}z^{2n},
 \end{equation}
where $$q_{n}=\frac{(-1)^{n}(\gamma)_{n+1}\Gamma(\beta)(2n+3)^{2}}{n!\Gamma(\alpha (n+1)+\beta)}\ \ \ \mbox{and}\ \ \ r_{n}=\frac{(-1)^{n}(\gamma)_{n}\Gamma(\beta)(2n+1)^{2}}{n!\Gamma(\alpha n+\beta)}.$$
By comparing  the coefficients of \eqref{rad19} and \eqref{rad20}  we obtain
$$\mu_{1}=\frac{9\gamma\Gamma(\beta)}{\Gamma(\alpha+\beta)},\quad \mu_{2}=\frac{81\gamma^2\Gamma^{2}(\beta)}{\Gamma^{2}(\alpha+\beta)}-\frac{25\gamma(\gamma+1)\Gamma(\beta)}{\Gamma(2\alpha+\beta)}$$
and by using the Euler-Rayleigh inequalities $$\mu_{k}^{-1/k}< \tau^2_{\alpha,\beta,\gamma,1}<\frac{\mu_{k}}{\mu_{k+1}}$$
for $k=1$ we have the inequalities of this part of the theorem.

{\bf b.} In view of the definition of the generalized Mittag-Leffler function we have
$$\varpi(\alpha,\beta,\gamma,z)=(zh^{\prime}_{\alpha,
\beta,\gamma}(z))^{^\prime}=1+\sum_{n \geq 1}\frac{(-1)^{n}\Gamma(\beta)(\gamma)_{n}(n+1)^{2}z^{n}}{n!\Gamma(\alpha n+\beta)}$$
and consequently
\begin{equation}\label{rad23}
\frac{\varpi'(\alpha,\beta,\gamma,z)}{\varpi(\alpha,\beta,\gamma,z)}=-\left.\sum_{n\geq 0}t_{n}z^{n}\right/\sum_{n\geq0}s_{n}z^{n},
\end{equation}
where
$$t_{n}=\frac{(-1)^{n}\Gamma(\beta)(\gamma)_{n+1}(n+2)^{2}}{n!\Gamma(\alpha (n+1)+\beta)}\ \ \
\mbox{and}\ \ \ s_{n}=\frac{(-1)^{n}\Gamma(\beta)(\gamma)_{n}(n+1)^{2}}{n!\Gamma(\alpha n+\beta)}.$$

Since $h_{\alpha,\beta,\gamma}\in \mathcal{LP},$ it follows that $h'_{\alpha,\beta,\gamma}\in\mathcal{LP},$ and consequently $z\mapsto \varpi(\alpha,\beta,\gamma,z)$ belongs also to the Laguerre-P\'{o}lya class $\mathcal{LP},$ and hence all its zeros are real. If we suppose that $\varsigma_{\alpha,\beta,\gamma,n}$ is the $n$th positive zero of the function $z\mapsto \varpi(\alpha,\beta,\gamma,z),$ then we get
$$\varpi\alpha,\beta,\gamma,z)=\prod_{n \geq 1}\left(1-\frac{z}{\varsigma_{\alpha,\beta,\gamma,n}}\right)$$
and for $|z|<\varsigma_{\alpha,\beta,\gamma,1}$
\begin{equation}\label{rad22}
\frac{\varpi'(\alpha,\beta,\gamma,z)}{\varpi(\alpha,\beta,\gamma,z)}=-\sum_{n\geq1}\frac{1}{\varsigma_{\alpha,\beta,\gamma,n}-z}
=-\sum_{n\geq1}\sum_{k\geq0}\frac{z^k}{\varsigma^{k+1}_{\alpha,\beta,\gamma,n}}=-\sum_{k\geq0}\sum_{n\geq1}\frac{z^k}{\varsigma^{k+1}_{\alpha,\beta,\gamma,n}}
=-\sum_{k \geq 0} \nu_{k+1}z^{k},
\end{equation}
where $\nu_{k}=\sum_{n \geq 1}\varsigma^{-k}_{\alpha,\beta,\gamma,n}.$
By comparing the coefficients of \eqref{rad22} and \eqref{rad23}, we have
$$\nu_{1}=\frac{4\gamma\Gamma(\beta)}{\Gamma(\alpha+\beta)}, \ \ \nu_{2}=\frac{16\gamma^2\Gamma^{2}(\beta)}{\Gamma^{2}(\alpha+\beta)}-\frac{9\gamma(\gamma+1)\Gamma(\beta)}{\Gamma(2\alpha+\beta)}$$
and by using the Euler-Rayleigh inequalities $$\nu_{k}^{-1/k}< \varsigma_{\alpha,\beta,\gamma,1}<\frac{\nu_{k}}{\nu_{k+1}}$$ for $k=1$ we have the inequalities of the second part of the theorem.
\end{proof}

\section*{Acknowledgements} The research of A. Prajapati was supported under the INSPIRE fellowship, Department of Science and
Technology, New Delhi, Government of India, Sanction Letter No. REL1/2016/2/2015-16.

\end{document}